\documentclass[11pt,leqno]{amsart}
\usepackage{graphicx}

\numberwithin{equation}{section}
\newcommand{\beqn}{\begin{eqnarray}}
\newcommand{\eeqn}{\end{eqnarray}}

\newcommand{\beq}{\begin{equation}}
\newcommand{\eeq}{\end{equation}}

\newcommand{ \db}{\bar{\partial}}

\newtheorem{theorem}{Theorem}[section]

\theoremstyle{definition}
\newtheorem{definition}[theorem]{Definition}

\theoremstyle{definition}

\theoremstyle{definition}

\theoremstyle{remark}

\begin{document}

\title{The Modified Calabi-Yau problems for CR-manifolds }

\author{Jianguo Cao}

\address{Mathematics Department, University of Notre Dame,
Notre Dame, IN 46556, USA; }

\address{Department of Mathematics, Nanjing
University, Nanjing 210093, China.}

\email{cao.7@nd.edu}

\author{ Shu-Cheng Chang}

\address{Department of Mathematics, National Tsing Hua University,
Hsinchu 30013, Taiwan, R. O. C.}

\email{scchang@math.nthu.edu.tw}


\date{}

\dedicatory{Dedicated to the memory of Xiao-Song Lin.}


\begin{abstract} In this paper, we derive a partial result
related to a question of Professor Yau: ``Does a simply-connected complete
K\"ahler manifold M with negative sectional curvature admit a
bounded non-constant holomorphic function?"

\medskip

\noindent
\textbf{Main Theorem.}
{\itshape
 Let $M^{2n}$ be a simply-connected complete
K\"ahler manifold M with negative sectional curvature $ \le -1 $
and $S_\infty(M)$ be the sphere at infinity of $M$. Then there is
an explicit {\it bounded } contact form $\beta$ defined on the
entire manifold $M^{2n}$.

Consequently, if $M^{2n}$ is a simply-connected K\"ahler manifold
with negative sectional curvature $-a^2 \le sec_{M} \le -1$, then
the sphere $S_\infty(M)$ at infinity of $M $ admits a {\it
bounded} contact structure and a bounded pseudo-Hermitian metric
in the sense of Tanaka-Webster.
}

 We also discuss several open modified problems of Calabi and
   Yau for Alexandrov spaces and CR-manifolds.
\end{abstract}
\maketitle

\section*{0. Introduction}

In this paper, we will provide a detailed construction of {\it bounded} contact structures on  a simply-connected complete K\"ahler manifold M
with negative sectional curvature $ \le -1 $. Afterwards, we will discuss the related open problems inspired by  Calabi and Yau.

 In 1979, Professor S. T. Yau [Y1] asked the following question.

 \bigskip

\noindent \textbf{Problem 0.1.} {\itshape (Yau [Y1])  Let $M^{2n}$ be a simply-connected complete K\"ahler manifold M with negative sectional curvature $
\le -1 $. Does $M^{2n}$ admit a bounded non-constant holomorphic function? }

\medskip

 In fact, an even more attractive problem in complex analytic differential geometry is to characterize bounded domains in $C^n$ within
 noncompact manifolds.

 \bigskip

\noindent
\textbf{Problem 0.2.}
{\itshape (Yau [Y1]) Let $M^{2n}$ be a simply-connected
complete K\"ahler manifold M with negative sectional curvature $ \le -1 $. Is $M$ bi-homeomorphic to a bounded domain in $\Bbb C^n$?
}

\medskip

  There were  some partial progress
  made  by Bland [Bl] and Nakano-Ohsawa [NO].
   Under extra assumptions, they  proved the existence of CR functions on the ideal boundary $S_\infty(M)$. In [Bl], two sufficient conditions were
    given for a complete K\"ahler manifold $M$ of non-positive sectional curvature to
  admit nonconstant bounded holomorphic functions, which seems also to guarantee that $M$ is a relatively compact domain with smooth boundary.

The precise definition of ideal boundary $S_\infty(M)$ can be
found in [BGS].

\bigskip

\noindent
\textbf{Theorem 0.3.}
{\itshape Let $M^{2n}$ be a simply-connected complete K\"ahler manifold M with negative sectional curvature $ \le -1 $ and
$S_\infty(M)$ be the sphere at infinity of $M$. Then there is an explicit {\it bounded } contact form $\beta$ defined on the entire manifold
$M^{2n}$.

Consequently, if $M^{2n}$ is a simply-connected K\"ahler manifold
with negative sectional curvature $-a^2 \le sec_{M} \le -1$, then
the sphere $S_\infty(M)$ at infinity of M admits a {\it bounded}  contact structure and a bounded pseudo-Hermitian metric in the
sense of Tanaka-Webster.
}

\medskip

Our proof of Theorem 0.3 was inspired by Gromov's bounded cohomology [Gro1-2] and calculations in [CaX].

Let $\omega $ be the K\"ahler metric on $M^{2n}$. It is clear that $d\omega = 0$. When $M^{2n}$  is a simply-connected complete K\"ahler
manifold with negative sectional curvature $ \le -1 $, Gromov observed that there must be a bounded $1$-form $\beta$ with
$$
   d\beta = \omega. \eqno(0.1)
$$

The proof of Gromov's assertion was outlined in [Pa] and [JZ]. In this paper, we provide a detailed proof of Gromov's assertion in \S 1. A
similar sub-linear estimates for equation (0.1) on manifolds with non-positive curvature was given by the first author and Xavier in [CaX].
\bigskip

\section*{1.   Bounded solutions to $d\beta=\alpha$ on manifolds with negative curvature }

In this section, we prove Theorem 0.3. In addition, we present a new direct proof of Gromov's bounded cohomology
theorem of negative curvature, see Theorem 1.4 and its proof below. Gromov's original approach to Theorem 1.4 below was based
a volume estimate of $k$-dimensional cone over a $(k-1)$-dimensional chain, and then use a dual space argument to complete
the proof. Our new method is to work on $k$-chains directly with a controlled Poincar\'e lemma for negative curvature.
Our approach might have some potential independent applications.

Throughout this section $(M^m, g)$ will be a complete simply-connected
manifold of negative sectional
curvature $ \le -1$. Let also $\alpha$  be a bounded smooth closed $k$-form on
$M$ with $k \ge 1$.
Since $M^m$ is diffeomorphic to
$\Bbb R^m$ there exists a form $\beta$ such that $d\beta =\alpha$.
The purpose of this section is to
 show that $\beta$ can be chosen to be bounded. The proof will
follow from the
Poincar\'e lemma by a comparison argument.

Fix $p\in M$ and denote by $exp_p: T_pM \rightarrow M$ the exponential map based
at $p$.

\bigskip

\noindent
\textbf{Lemma 1.1.}
{\itshape
  Consider the maps $\tau_t:M\rightarrow M$, given by
$x \longmapsto exp_p(t \,  exp_p^{-1}(x))$, where $0\leq t\leq 1$. Then
$$
|(\tau_t)_{\ast}\xi|\leq \frac{\sinh tr}{\sinh r} |\xi| \eqno(1.1)
$$
for every tangent vector $\xi$, where $r = d(x, p)$.}
\begin{proof} Let $\sigma: [0, 1] \rightarrow M^n$ be the geodesic segment joining
$p$ to $x$,
$ \xi \in T_xM^n$ and $y = (exp_p)^{-1}(x) \in   T_pM^n$.
By a straightforward computation one has
$$
\aligned
  & (\tau_t)_* \xi
= (d\; exp_p)_{t (exp_p)^{-1}(x)} [t d (exp^{-1}_p)_{(x)} \xi ] \\
& = (d \; exp_p)_{ty} \{t [d (exp_p)_y]^{-1} \xi\}.
\endaligned
$$

  Recall that $\sigma(t) = exp_p(ty)$.  It is now manifest from the above formula that
$$
     J(tr) := (\tau_t)_* \xi  \eqno(1.2)
$$
 is the Jacobi field along $\sigma$  satisfying
$J(0)=0,\;J(r)=  \xi$. On the other hand, since the sectional curvatures are $\le -1$,  we estimate the function $f(s):=|J(s)|$ by a method inspired by
Gromov. It is sufficient to verify
$$
     \frac{|J(s)|}{\sinh s} \le   \frac{|J(r)|}{\sinh r}, \eqno(1.3)
$$
for all $0 \le s \le r$.
\medskip

We may assume that $r > 0$, otherwise the inequality (1.1) holds trivially.
To do this, we consider the function
$$
   \eta(s) = \frac{f(s)}{\sinh s}.
$$

It is sufficient to verify
$$
\frac{f(s)}{\sinh s} \le \frac{f(r)}{\sinh r} \text{  or   } \eta'(s) \ge 0. \eqno(1.4)
$$
Since we have
$$
\eta'(s) = \frac{f'(s) \sinh s - f(s) \cosh s }{[\sinh s]^2},
$$
it remains to verify that
$$
[f'(s) \sinh s - f(s) \cosh s]' = f''(s) \sinh s - f(s) \sinh s \ge 0. \eqno(1.5)
$$

Recall that the curvature tensor $R$ is given by $R(X, Y) Z = - \nabla_X \nabla_Y Z + \nabla_Y \nabla_X Z + \nabla_{[X, Y]}Z$
where $[X, Y] = XY - YX$ is the Lie bracket of $X$ and $Y$.

Following a calculation in  [BGS],  by our assumption of $sec_{M} \le -1$ we have
$$
\aligned
f''(s) & = |J(s)|'' \\
&  = [ \frac{\langle J, J' \rangle }{|J|}]' \\
& = \frac{ \langle J, J'' \rangle |J|^2 + \langle J', J' \rangle |J|^2 -  \langle J, J' \rangle^2  }{|J|^3} \\
& \ge  \frac{- \langle R(\sigma', J)\sigma', J \rangle |J|^2   }{|J|^3} \\
& \ge f(s),
\endaligned \eqno(1.6)
$$
where we used the assumption that $ \langle J'', J \rangle = - \langle R(\sigma', J)\sigma', J \rangle \ge |J|^2$.
It follows from (1.5)-(1.6) that (1.4) holds. This completes the proof of (1.3) as well as Lemma 1.1.
\end{proof}

Recall that if $\alpha$ is a
$k$-form and $Z$ is a vector field,  then  $(\alpha\lfloor_{ Z})$
is the $(k-1)$-form given by
$$(\alpha\lfloor_{ Z})(\xi_1,\cdots,\xi_{k-1}) = \alpha(
Z, \xi_1,\cdots,\xi_{k-1}). $$

 For the sake of completeness we give a proof of the following elementary result.

\bigskip

\noindent \textbf{Lemma 1.2.} {\itshape
 Let $\Psi$ be a closed $k$-form in $\Bbb R^m$.
Then the $(k-1)$-form $\Phi$ defined by
$$\Phi(x) = r \int^1_0
[(\tau_t)^*(\Psi\lfloor_{\frac{\partial}{\partial r}})](x)dt
$$
satisfies $d\Phi=\Psi$; here
${\frac{\partial}{\partial r}} =
 \sum^m_{i=1} \frac{x_i}{r}  \frac{\partial}{\partial x_i}$, $r =
(\sum^m_{i=1}x^2_i)^{1/2}$
 and
$\tau_t(x)=tx$.}
\begin{proof} By the standard proof of the Poincar\'e lemma
([SiT], p.130), $\Phi$ can be taken to
be
$$
\Phi(x)=\sum\limits_{i_1<\cdots<i_k}
\sum\limits^k\limits_{j=1} (-1)^{j-1} x_{i_j}
\Big( \int^1_0 t^{k-1}
\Psi_{i_1\cdots i_k} (tx) dt \Big) dx_{i_1}\wedge\cdots\wedge\widehat{dx_{i_j}}
\wedge\cdots\wedge dx_{i_k},
$$
where $\Psi =
\sum_{i_1<\cdots<i_k}  \Psi_{i_1\cdots
i_k}  dx_{i_1}\wedge\cdots\wedge dx_{i_k}$.

In particular, one has
$$
\aligned
 \Phi (x) & =  \sum\limits_{i_1<\cdots<i_k}
\sum\limits^k\limits_{j=1}  x_{i_j}
\Big( \int^1_0 t^{k-1}
\Psi_{i_1\cdots i_k} (tx) dt \Big) (dx_{i_1}\wedge\cdots\wedge dx_{i_k})
\lfloor_{\frac{\partial}{\partial x_{i_j}}} \\
& = r \sum\limits_{i_1<\cdots<i_k}\Big( \int^1_0 t^{k-1}
\Psi_{i_1\cdots i_k} (tx) dt \Big) (dx_{i_1}\wedge\cdots\wedge dx_{i_k})
\lfloor_{\frac{\partial}{\partial r}} \\
& = r \int^1_0 t^{k-1} (\Psi\lfloor_{\frac{\partial}{\partial r}})(tx)dt \\
& = r \int^1_0
[(\tau_t)^*(\Psi\lfloor_{\frac{\partial}{\partial r}})](x)dt,
\endaligned
$$
as desired. \end{proof}

   We would also like to borrow another elementary but useful observation of Gromov, in order to prove our main theorem

\bigskip

\noindent
\textbf{Lemma 1.3.}
{\itshape
(Gromov, [Cha, page 124]) Suppose that $f$ and $h$ are positive integrable functions, of real variable $r$,
for which
$$
\frac fg
$$
is an increasing with respect to $r$. Then the function
$$
\frac{\int_0^r f }{\int_0^r g}
$$
is also increasing with respect to $r \ge 0$.}

Let us now provide a new detailed proof of a theorem of Gromov.

\bigskip

\noindent
\textbf{Theorem 1.4.}
{\itshape (Gromov) Let $M^m$ be a simply-connected complete Riemannian manifold with negative sectional curvature $\le -1 $.
Suppose that $ \alpha$ is bounded closed $k$-form with $k\ge 2$. There is a bounded $(k-1)$-form $\beta$ with $
   d \beta = \alpha
$
satisfying
$$
   \| \beta \|_{L^\infty} \le \frac{1}{k-1}    \| \alpha \|_{L^\infty}. \eqno(1.7)
$$}
\begin{proof}
Let $(x_1,..., x_n)$ be Euclidean coordinates in  $T_pM$ and
consider the pull-back metric $h$ of the metric $g$ under $exp_p:T_pM \rightarrow M$.
Observe that there are now two ways to interpret the map $\tau_t$.
The first interpretation comes from Lemma 1.1 with $(M, g)$ being replaced by
$(T_pM, h)$;
alternatively, one can think of $\tau_t$ as
the self-map of $T_pM$, $(x_1,..., x_n)\longmapsto t(x_1,..., x_n)$,
that appears in the Poincar\'e lemma (Lemma 1.2).
It is an  easy and  yet basic observation that
these two ways of thinking about
$\tau_t$  give rise to the same map.

  We may also replace the form $\alpha$ that appears in the statement of
Lemma 1.2 by a closed form $\Psi$ on $T_pM$ which is bounded in
the induced  metric $h$. Let $\Phi$ be given by Lemma 1.2 and observe that,
by Lemma 1.1,
$$
    |(\tau_t)^* \varphi (x) |_h \leq \big( \frac{\sinh tr }{\sinh r } \big)^{k-1} |\varphi(\tau_t(x))|_h, \; \; \;
  k \geq 2,  \eqno(1.8)
$$
holds for any $(k-1)$-form $\varphi$ on $T_pM$; here $|\cdot|_h$ is any
one of the equivalent  norms induced by $h$.
 Since $|\frac {\partial}{\partial r}| = 1$, it follows from (1.3) and Lemma 1.2
that
$$
\aligned  | \Phi (x)|_h
& \leq r \int^1_0 |[(\tau_t)^*(\Psi\lfloor_{\frac{\partial}{\partial r}})](x)|_h dt \\
& \leq
 r \int^1_0 \big( \frac{\sinh tr }{\sinh r } \big)^{k-1}  | \Psi  (tx)\lfloor_{\frac{\partial}{\partial r}} |_h dt
\\
& = \int^r_0 \big( \frac{\sinh s }{\sinh r } \big)^{k-1}  | \Psi
(\frac sr x)\lfloor_{\frac{\partial}{\partial r}} |_h ds
\\
&\leq \frac{\int^r_0 (\sinh s)^{k-1} ds }{(\sinh r)^{k-1} }
  \sup_{0 \le s \leq r}
 |\Psi( \frac sr x) |_h
 \endaligned \eqno(1.9)
$$

Choosing $f(r) = (\sinh r)^{k-1}$ and $\hat g(r) = (k-1) (\sinh r)^{k-2} \cosh r$ in Lemma 1.3, we see that $[\frac{f}{\hat g}]' =
\frac{1}{(k-1) (\sinh r)^2} > 0$ and
$$
\frac{\int^r_0 (\sinh s)^{k-1} ds }{(\sinh r)^{k-1} } \le \frac{1}{k-1}. \eqno(1.10)
$$
It follows from (1.9)-(1.10) that
$$
  |\Phi(x)|_h \leq  \frac{1}{k-1}  \sup |\Psi|_h. \eqno(1.11)
$$
Hence  $\Phi$ is a bounded solution of $d  \Phi  = \Psi$ and the proof of Theorem 1.4 is completed.  \end{proof}

\bigskip
\noindent
\textbf{Proof of Main Theorem:}

Our main theorem (Theorem 0.3.) can be derived as follows. We fix a base point $p$ as above. There is a differential structure $\Xi_p$
imposed on $S_\infty(M)$ given by the map
$$
\aligned
F_p  : & \overline{ B_1(0)} \to M \cup S_\infty(M) \\
&  \vec{v} \to Exp_p [ \frac{\vec{v}}{1 - | \vec{v}| }      ].
\endaligned
$$

For $p \neq q$, the transitive map $F_q^{-1} \circ F_p: \overline{ B_1(0)} \to \overline{ B_1(0)}$ is not necessarily smooth.
However, we fix {\it one} differential structure $\Xi_p$ on $S_\infty(M)$ via the map $F_p$.

Let $J$ be the complex structure of our K\"ahler manifold $M$. Let $r(x) = d(x, p)$ and $\beta =  J \circ dr$, i.e.,
$\beta(\vec w) = dr ( J \vec w) $ for all $\vec w \in T_x(M)$.  When $- a^2 \le sec_M -1 $, it is known that
$$
    |X|^2 \le  | (\nabla_X dr) (X)| = |Hess(r) (X, X)|  \le  a |X|^2
$$
for all $X \in T_x(\partial B_r(p))$ with $r >> 1$.

  Since $M$ is K\"ahler, we have $\nabla_X J = 0$. It follows that $|\nabla_X \beta| \le a |X|$ for $X \in T_x(\partial B_r(p))$ with $r >> 1$.

  Thus, $\{ \beta |_{\partial B_r(p)} \}$ defines an equi-continuous family of contact forms on $S_\infty(M)$. By Ascoli lemma,
   there
  is a subsequence to converge to a bounded contact form $\beta_\infty$ on
  $S_\infty(M)$. Since $sec_M \le -1$, it is known that $ d\beta (\tilde X, \bar {\tilde X}) = Hess(r) (X,
  X) + Hess(r) (JX, JX) \ge 2|X|^2 $ for all $X \in T_x(\partial B_r(p))$ and $X \perp \nabla r$, where
  $\tilde X = \frac{1}{\sqrt{2}}[ X - \sqrt{- 1} JX ]$. Therefore, $\beta_\infty$ defines a non-trivial contact form on $S_\infty(M)$.
  Moreover, $\omega_\infty = d \beta_\infty$ gives rise to a pseudo-hermitian metric on $S_\infty(M)$.

     Similarly, one can also choose $\beta^*$ satisfying  $d \beta^* = \omega$, where $\omega$ is the K\"ahler form of $M$ and $\beta^*$ in the proof
     of Theorem 1.4. With extra efforts, one can show that $|\nabla \beta^* | \le c_1$ for some constant $c_1$.
     Thus, $\{ \beta^* |_{\partial B_r(p)} \}$ defines an equi-continuous family of contact forms on $S_\infty(M)$ as well.

     This completes the proof of our main theorem.

\section*{2. The modified Calabi-Yau problems for singular spaces and CR-manifolds}

In this section, we will discuss the generalized Calabi problems on K\"ahler
manifolds with boundaries. In addition, we will comment on the existence of positive sup-harmonic functions on
(possibly singular) Alexandrov spaces with non-negative sectional
curvature.

\bigskip

\noindent
\textbf{\S A. Sup-harmonic functions on Alexandrov spaces with
non-negative sectional curvature }

\bigskip

Professor S. T. Yau also had earlier results on bounded harmonic functions on smooth complete Riemannian
manifolds with non-negative Ricci curvature. We would like to extend this theorem of Yau to singular spaces.

In an important paper [Per1], Perelman provided an affirmative
solution to the Cheeger-Gromoll soul conjecture. More precisely,
he showed that {\it ``if a smooth complete non-compact Riemannian
manifold $M^n$ of non-negative curvature has a point $p_0$ with
strictly positive curvature $K|_{p_0} > 0$, then $M^n$ must be
diffeomorphic to $\Bbb R^n$}. In [Per1], Perelman also asked to
what extent the conclusions of his paper [Per1] would hold for
(possibly singular) Alexandrov spaces with non-negative curvature.

Recently, the first author together with Dai and Mei showed the following.

\bigskip

\noindent
\textbf{Theorem A.1.}
{\itshape (Cao-Dai-Mei, 2007, [CaMD1]) Let
$M^n$ be a $n$-dimensional complete, non-compact Alexandrov space
with non-negative sectional curvature. Suppose that $M^n$ has no
boundary and $M^n$ has positive sectional curvature on an
non-empty open set. Then $M^n$ is contractible.
}

In 1976, Professor S. T. Yau proved the following Liouville type
theorem.

\bigskip

\noindent \textbf{Theorem A.2.} {\itshape (Yau, 1976 [Y3]) Let
$M^n$ be a $n$-dimensional complete, non-compact smooth Riemannian
space
 with non-negative Ricci curvature.  Then any positive harmonic functions on $M^n$  must be a constant function.}

On an (possibly singular) Alexandrov space, we introduce the
following notion of sup-harmonic function.

\definition{Definition A.3} Let $M^n$ be a $n$-dimensional complete, non-compact Alexandrov
space with non-negative sectional curvature. Suppose that $M^n$
has no boundary,  $f: M^n \to \Bbb R$ is a Lipschitz continuous
function and
$$
    f(x) \ge \frac{1}{Area(\partial B_\varepsilon(x))} \int_{\partial B_\varepsilon(x)}  f dA \eqno(A.1)
$$
for any sufficiently small $\varepsilon >0$. Then we say that $f$
is a sup-harmonic function on $M$.
\enddefinition

For example, $ f(x) = - [d(x, x_0)]^2$ is a sup-harmonic function
on $M$, whenever $M$ has non-negative sectional curvature in
generalized sense.

\bigskip

\noindent \textbf{Problem A.4.} {\itshape (Liouville-Yau type
problem) Let $M^n$ be a $n$-dimensional complete, non-compact
Alexandrov space with non-negative sectional curvature. Suppose
that $M^n$ has no boundary.   Is it true that  any positive
sup-harmonic functions on $M^n$  must be a constant function? }

In [CaB], the first author and Benjamini  studied a different
Liouville-type problem of Schoen-Yau. One hopes to continue to
work on Liouville-Yau type problem mentioned above.

\bigskip

\noindent
\textbf{\S B. The generalized Calabi problems for K\"ahler
domains with boundaries}
\medskip

The classical Calabi problems for Ricci curvatures on  compact
K\"ahler manifolds {\it without boundaries} have been successfully
solved by Professor S. T. Yau.

\bigskip

\noindent
\textbf{Theorem B.1.}
{\itshape
(Yau [Y4]) Let $M^{2n}$ be a compact smooth
K\"ahler manifold without boundary. Then the following is true:
\smallskip
\noindent (1) For any K\"ahler form $\omega_0 \in H^{(1,
1)}(M^{2n})$ and any $(1, 1)$-form  $\beta$ representing the first
Chern class $c_1(M^{2n})$, there is a K\"ahler metric $\tilde
\omega = \omega_0 + i \partial\db f $ such that its Ricci
curvature tensor satisfies
$$
Ric_{\tilde \omega}= \beta;
$$

\smallskip
\noindent (2) If the first Chern class $c_1(M) \le 0$, then
$M^{2n}$ admits a K\"ahler-Einstein metric.
}

For a K\"ahler manifold $\Omega$ with boundary $M^{2n-1} =
b\Omega$, we consider a similar problem. This problem is closely
related to the existence problem of CR-Einstein metrics, or
partially Einstein metrics.

\bigskip

\noindent
\textbf{Definition B.2.}
{\itshape
 (CR-Einstein metrics or partially Einstein metrics, [Lee2]) Let $\Sigma^{2n-1}$ be a $CR$-hypersurface
with $CR$-distribution $ \mathcal H_{\Sigma^{2n-1}} = \ker \theta$ for
some contact $1$-form $\theta$ and let $g_\theta (X, JY) =
d\theta (X, JY)$ be a pseudo-hermitian metric as above. If the
Ricci tensor of $g_\theta$ satisfies
$$
   Ric_{g_\theta} (X, Y) = c g_\theta (X, Y)
$$
for all $X, Y \in \mathcal H_{\Sigma^{2n-1}} = \ker \theta$ where $c$
is a constant, then $g_\theta$ is called a CR-Einstein (partially
Einstein) metric.
}

Inspired by Yau's result, Lee proposed to study the $CR$-version
of the Calabi problem.

\bigskip

\noindent
\textbf{Problem B.3.}
{\itshape
($CR$-Calabi Problems, [Lee2]) Let
$M^{2n-1}$ be a $CR$-manifold,  $\Phi$  be a closed form
representing the first Chern class for the bundle $T^{(1,
0)}(M^{2n-1})$ and $\Phi_b(X, Y) = \Phi (X, Y)$ for $X, Y \in \mathcal
H_{\Sigma^{2n-1}} = \ker \theta$.

 (1) Can we
find a pseudo-metric $g_{ \theta}$ such that its Ricci tensor
satisfies
$$
   Ric_{g_\theta} (X, Y) = \Phi_b (X, Y) \eqno(B.1)
$$
for all $X, Y \in \mathcal H_{\Sigma^{2n-1}} = \ker \theta$?

(2) Given a $(1, 1)$-form $\beta_b \in [c_1(M^{2n-2}]_b$,  can we
find a pseudo-metric $g_{ \theta}$ such that its Ricci tensor
satisfies
$$
   Ric_{g_\theta} (X, Y) = \beta(X, Y) \eqno(B.2)
$$
for all $X, Y \in \mathcal H_{\Sigma^{2n-1}} = \ker \theta$?
}

The pseudo-Hermitian metric for general $CR$-manifolds was also
discussed  in [Ta1-2] and [Web]. Authors  derived the
following partial answer to  Problem 3:

\bigskip

\noindent
\textbf{Problem B.4.}
{\itshape ([CaCh]) Let $M^{2n-1}$ be the smooth
boundary of a bounded strongly pseudo-convex domain $\Omega$ in a
complete Stein manifold $V^{2n}$.   Then for $n \ge 3$, $M^{2n-1}$
admits a CR-Einstein metric (or partially Einstein metric).
}

One might be able to continue working  on
Problem B.3, using Kohn-Rossi's $\db_b$-theory described below.

\bigskip

\noindent
\textbf{\S C. The Calabi-Escobar type problem for K\"ahler
domains with boundaries}
\medskip

The first author and Mei-Chi Shaw studied the $CR$-version of the
Poincar\'e-Lelong equation $i \partial_b \db_b u = \Psi_b$ in
[CaS3]. The linearization equation for (B.2) is related to the
$CR$-version of Poincar\'e-Lelong equation.

   In fact, to solve the linear function
   $$
\db_b u = \beta_b \text{  on   } b\Omega, \eqno(C.1)
   $$
Kohn and Rossi [KoRo] used the solutions to the $\db$-Cauchy
problem to solve $ \db_b u = \beta_b$ extrinsically as follows.
Let us first choose an arbitrary smooth extension $\hat \beta$ on
$\Omega$. If we can solve
$$
\begin{cases}\aligned
& \db v = \db \hat \beta  \text{  on   } \Omega \\
&  v\llcorner_X = 0,  \text{  for    } X \in T_z^{(0, 1)}(b\Omega)
\endaligned
\end{cases} \eqno(C.2)
$$
Clearly $\tilde \beta = \hat \beta - v $ is a $\db$-closed
extension on $\Omega$ of $\beta$. If we solve
$$
\db \tilde u  =  \hat \beta - v   \text{  on   } (\Omega \cup
b\Omega) , \eqno(C.3)
$$
then the restriction $u = \tilde u \llcorner_{b\Omega  }$
satisfies
$$
 \db_b [(\tilde u) \lfloor_{b\Omega  }]  = \beta_b \text{  on   } b\Omega.
$$

The details for solving the $\db$-Cauchy problem (C.2) was given
in Chapter 9 of [ChSh].

   In 1992, Escobar [Esc] was able to solve the non-linear curvature equation on manifolds {\it with boundary}.

\bigskip

\noindent
\textbf{Theorem C.1.}
{\itshape (Escobar [Esc]) Let $\Omega \subset \Bbb
R^n$ be a compact domain with smooth boundary $\partial \Omega$
and $n > 6$. Then there is a conformally flat metric $g$ on
$\Omega$ such that the scalar curvature $Scal_g$ of $g$ is zero
and the mean curvature $H_g$ of $(\partial \Omega, g)$ is
constant:
$$
\begin{cases} \aligned
& Scal_g = 0 \text{  on   } \Omega \\
&  H_g  = c \text{  on   } \partial\Omega,
\endaligned
\end{cases} \eqno(C.4)
$$
for some constant $c$.
}

Inspired by Theorem C.1 and the Kohn-Rossi's solution to
$\db$-Cauchy problem, we are interested in the following type.

\bigskip

\noindent
\textbf{Problem C.2.}
{\itshape (Calabi-Escobar type problem) Let $\Omega$
be a compact domain in Stein manifold $M$ with smooth strongly
pseudo convex boundary $b\Omega$, and let $H^{CR}_g$ be the
partial sum of second fundamental form of $(b\Omega, g)$ over the
CR-distribution $ker \theta$ of $b\Omega$. Is there is an
K\"ahler-Einstein metric $g $ on $\Omega$ with constant CR-mean
curvature on the boundary $b\Omega$? In another words, we would
like to find the existence of solution to the following non-linear
boundary problem:
$$
\begin{cases} \aligned
& Ric_g = c_1 g \text{  on   } \Omega \\
&  H^{CR}_g  = c_2 \text{  on   } b\Omega
\endaligned
\end{cases} \eqno(C.5)
$$
for some constant numbers $c_1$ and $c_2$.
}

The linearalization of non-linear equation is the Poincare-Lelong
equation with boundary conditions. The first author and Mei-Chi Shaw [CaS]
were able to solve
$$
i \partial_b  \db_b u = \Theta_b \mbox{  \quad on  \quad } b\Omega \eqno(C.6)
$$
even for weakly pseudo-convex domains $\Omega$ in $\Bbb CP^n$.

One hopes to
continue to work in direction, in order to investigate Problem
C.2.

\subsection*{Acknowledgment}

The first author is supported in part by an NSF grant Grant
DMS0706513 and  Changjiang Scholarship of China at Nanjing
University.  The first named author is also  very grateful to
National Center for Theoretical Sciences at National Tsinghua
University for its warm hospitality. The second author is
supported partially by the NSC of Taiwan.

   We thank Professor Mei-Chi Shaw for pointing a misquote in the earlier version of
   Section 1. She also proposed to use (not necessarily unique) one of smooth differential structures at infinity. The first
   author is indebted to Professor Rick Schoen for his advices on open manifolds of negative curvature.


\end{document}